\documentclass[]{article}

\usepackage[utf8]{inputenc}
\usepackage[margin=1in]{geometry}
\usepackage[numbers]{natbib}

\usepackage{suffix}
\usepackage{graphicx}
\usepackage[hidelinks]{hyperref}
\usepackage{setspace}
\usepackage{amsmath}
\usepackage{amssymb}
\usepackage{amsthm}
\usepackage{mathtools}
\usepackage{xcolor}
\usepackage{paralist}
\usepackage{caption}
\usepackage{subcaption}
\usepackage{booktabs}
\usepackage{enumitem}
\usepackage{setspace}
\usepackage[textsize=tiny,backgroundcolor=pink,prependcaption=TODO: ]{todonotes}
\allowdisplaybreaks
\newcommand{\RR}{\mathbb{R}}
\newcommand{\KK}{\mathbb{K}}
\newcommand{\NN}{\mathbb{N}}
\newcommand{\ZZ}{\mathbb{Z}}
\newcommand{\EE}{\mathbb{E}}

\newcommand{\XX}{\mathbb{X}}
\renewcommand{\AA}{\mathbb{A}}
\newcommand{\PP}{\mathbb{P}}

\newcommand{\one}[1]{\mathbf{1}_{\set{#1}}}
\WithSuffix\newcommand\one*[1]{\mathbf{1}_{#1}}

\renewcommand{\bar}{\overline}

\newcommand{\revision}[1]{#1}
\DeclarePairedDelimiter{\set}{\{}{\}}
\DeclarePairedDelimiter{\abs}{|}{|}

\DeclareMathOperator{\lev}{\mathcal{D}}
\DeclareMathOperator{\dom}{\mathrm{dom}}
\DeclareMathOperator{\Gr}{\mathrm{Gr}}
\DeclareMathOperator*{\argmin}{\mathrm{arg\,min}}
\newtheorem{theorem}{Theorem}[section]
\newtheorem{proposition}[theorem]{Proposition}
\newtheorem{lemma}[theorem]{Lemma}
\newtheorem{corollary}[theorem]{Corollary}

\theoremstyle{definition}
\newtheorem{definition}[theorem]{Definition}

\newtheorem{example}[theorem]{Example}

\theoremstyle{remark}
\newtheorem{remark}{Remark}

\newtheoremstyle{special}{}{}{\itshape}{}{\bfseries}{.}{.5em}{#1 (\thmnote{#3})}
\theoremstyle{special}
\newtheorem*{specialAssumption}{Assumption}

\makeatletter
\newcommand{\setword}[2]{%
	\phantomsection
	#1\def\@currentlabel{\unexpanded{#1}}\label{#2}%
}
\makeatother

\newlength\caseLen
\newlist{enumcase}{enumerate}{1}
\setlist[enumcase,1]{
	label=\textit{Case~\arabic*:},
	wide,
	nosep
}

\newcommand{\email}[1]{\href{mailto:#1}{\texttt{#1}}}

\title{Continuity of Discounted Values and the Structure of Optimal Policies for Periodic-Review Inventory Systems with Setup Costs}

\date{\today}
\author{
	Eugene A. Feinberg\thanks{Department of Applied Mathematics and Statistics, State University of New York at Stony Brook, Stony Brook, NY 11794-3600, USA, %
		\email{eugene.feinberg@stonybrook.edu}} \and
	David N. Kraemer\thanks{Department of Applied Mathematics and Statistics, State University of New York at Stony Brook, Stony Brook, NY 11794-3600, USA, %
		\email{david.kraemer@stonybrook.edu}}
}

\begin{document}
\maketitle
\begin{abstract}
	This paper proves continuity of value functions in discounted periodic-review single-commodity total-cost inventory control problems with \revision{continuous inventory levels,} fixed ordering costs, possibly bounded inventory storage capacity, and possibly bounded order sizes for finite and infinite horizons. In each of these constrained models, the finite and infinite-horizon value functions are continuous, there exist deterministic Markov optimal finite-horizon policies, and there exist stationary deterministic Markov optimal infinite-horizon policies. For models with bounded inventory storage and unbounded order sizes, this paper also characterizes the conditions under which $(s_t, S_t)$ policies are optimal in the finite horizon and an $(s,S)$ policy is optimal in the infinite horizon. \\

\end{abstract}
\section{Introduction}
\label{sec:introduction}

Periodic review inventory control studies have a rich history in operations research beginning with~\citet{arrow_1951_optimal} and~\citet{dvoretzky_1952_inventory}. The classical results, including on multistage problems, are summarized in the texts~\citet{bensoussan_2011_dynamic,porteus_2002_foundations,simchi_levi_2014_logic,zipkin_2000_foundations}. One of the principal methods of studying inventory control problems is the analysis of $(s,S)$ policies. Under an $(s,S)$ policy, the controller only orders when the inventory level drops below $s$, and the amount that is ordered returns the inventory level to $S$. Since the seminal paper of~\citet{scarf_1960_optimality}, $(s,S)$ policies have been shown to be optimal in a variety of problem formulations; see, e.g.,~\citet{veinott_1965_computing} for the case of discrete demand,~\citet{zabel_1962_note} and~\citet{iglehart_1963_optimality} for the case of continuous demand. Recent developments in the theory of Markov decision processes (MDPs) have improved the understanding of the periodic-review single-commodity inventory control problem with \revision{continuous inventory levels,} setup costs, unbounded storage capacity, unbounded order sizes, and backordering; see~\citet{feinberg_2016_optimality, feinberg_2018_convergence}. A complete description of the optimality of $(s,S)$ policies for convex holding costs in this setting is given in~\citet{feinberg_2017_structure}, where additional references can be found. For the optimality of $(s,S)$ policies under the average-cost criteria, see~\citet{feinberg_2018_convergence}.

This paper studies continuity of values under two additional constraints to the inventory problem with fixed order costs: possibly bounded storage capacity and possibly bounded order sizes. Bounded storage models were first considered in~\citet{hartley_1976_operations}, but remain an active field of research; see~\citet{jiang_2021_partial} and the references therein. Studies of models with bounded orders~(sometimes called capacitated models) include~\citet{shaoxiang_1996_x-y, xie_1998_stability,gallego_2000_capacitated, bensoussan_2007_partially,chao_dynamic_2012, bartoszewicz_2019_sliding}. Discount-optimal policies for models with bounded orders are not yet fully characterized. We also consider the version of the model with lost sales, in which unrealized demand is lost.~\citet{feinberg_2017_structure} proved continuity of the finite and discounted infinite horizon value functions for the setup cost inventory models with unbounded storage, unbounded order sizes, and backorders by studying optimal policies and in particular, optimal $(s,S)$ policies. One difficulty that models with possibly bounded storage and possibly bounded order sizes and models with lost sales present is that, for these models, $(s,S)$ policies may not be optimal or even feasible, and continuity of value functions must be established via other means.

Continuity of values is an important property for practical applications in two major ways. For approximate methods based on discretization of continuous state and action spaces, continuity of values is often required for convergence of solutions to the approximate problems to the solutions of the true problem as the discretization is refined. In addition, for problems where the controller only has access to a noisy or otherwise partially observable system transition model, continuity of values is necessary for empirical consistency of the solutions to the estimated models. For general overviews of these and related issues, see~\citet{kara_2020_robustness,kara_2021_qlearning}.

\revision{%
  The main results of this paper concern continuity of values in the infinite horizon versions of these problems. Theorem~\ref{thm:continuity_of_finite_horizon_values}, whose main results were obtained in~\cite[Theorem 14]{feinberg_2021_continuity}, establishes continuity of values in the finite horizon. Theorem~\ref{thm:continuity_of_infinite_horizon_values} proves continuity of values in the infinite horizon after showing an important inequality between the values for the inventory control problems with and without setup costs. In the finite horizon, a generalized form of Berge's maximum theorem is the main mechanism for establishing continuity of values; see~\citet[Theorems 8, 14]{feinberg_2021_continuity}. For the infinite horizon, continuity requires additional analysis to ensure that the the tail expectations of the finite values vanishes in the limit. These general facts apply independently of storage or order constraints. While the case of backorders with unbounded storage capacity was studied in \cite{feinberg_2017_structure}, the results for models with backorders with limited storage capacity and models with lost sales are new.  Theorem~\ref{thm:bounded_storage_s_S_policies} characterizes the optimality of $(s_t, S_t)$ and $(s,S)$ policies for problems with backorders and bounded storage capacity.
}

This paper is organized as follows. Section~\ref{sec:model_description} defines the stochastic periodic-review single-commodity inventory control problem with possibly bounded storage capacity and possibly bounded order sizes, \revision{and it provides results on continuity of values in the finite and the infinite horizon setting.} Section~\ref{sec:bounded_storage} considers inventory models with unbounded order sizes, and it characterizes the conditions under which $(s_t,S_t)$ and $(s,S)$ policies are optimal in the finite and infinite horizon, respectively.
\section{Model Description \revision{and Main Results}}
\label{sec:model_description}

\revision{
In this section we define the inventory models under consideration throughout the paper. We also prove the continuity of values in the finite and discounted infinite-horizon problems.
}
Let $\RR$ be the real numbers, $\RR_+ := [0, +\infty)$ the nonnegative real numbers, $\NN := \set{0, 1, 2, \dots}$ the natural numbers, \revision{and $\ZZ$ the integers}. The \emph{stochastic periodic-review setup-cost inventory model} is defined as follows. At times $t \in \NN$, a controller views the current inventory \revision{$x_t \in \RR$} of a single commodity and makes an ordering decision \revision{$a_t \geq 0$}. Lead times are zero, so the orders are filled instantaneously prior to the realization of demand. The cost of ordering is paid at the time of delivery of the order. For problems with backorders, any demand is satisfied and unment demand is backlogged. For problems with lost sales, the demand is satisfied up to the available inventory level, and if the demand is greater than this level, the unmet demand is lost. After the demand is satisfied, the controller views the remaining inventory and pays holding (i.e., excess inventory) or backorder (i.e., negative inventory) costs, and the process continues. The demand and order quantity are assumed to be nonnegative.

Let $a \land b := \min\{a,b\}$ and $a \lor b := \max\{a,b\}$. The inventory model is characterized by the following parameters:
\begin{enumerate}
	\item $K > 0$ is a fixed ordering cost, paid whenever the order size is nonzero;
	\item $\bar c > 0$ is the per-unit ordering cost;
	\item $h : \RR \to \RR_+$ is convex and continuous holding/backorder cost function with $h(x) \to +\infty$ as $\abs{x} \to +\infty$, where there is no loss in generality to assume that $\inf h = 0$;
	\item $\set{D_t : t = 1, 2, \dots}$ is a sequence of i.i.d. nonnegative finite random variables representing the demand at periods $1, 2, \dots$, where for such a random variable $D$, we assume that $P(D > 0) > 0$ and $\EE h(x-D) < +\infty$ for each $x \in \RR$;
	\item $T : \RR \to \RR$ is the backorder rule given by $T(x) := 0 \lor x$ for the inventory model with lost sales, and $T(x) := x$ for the inventory model with backorders;
	\item $\alpha \in [0, 1)$ is the discount factor;
	\item $\bar a > 0$ is the maximum order size, where $\bar a = +\infty$ means any finite order size may be placed;
	\item $\bar x > 0$ is the maximum storage size, where $\bar x = +\infty$ means any finite amount of inventory may be stored.
\end{enumerate}
The model with $\bar a = +\infty$ ($\bar a < +\infty$) is called an unbounded (bounded) orders model. The model with $\bar x = +\infty$ ($\bar x < +\infty$) is called an unbounded (bounded) storage model. A common example of a function $h(x)$ satisfying the above conditions is $h(x) = (0 \lor \eta_1 x) - (0 \land \eta_2 x)$, where $\eta_1, \eta_2 > 0$.

\begin{figure}[t]
	\centering
	\includegraphics[scale=0.25]{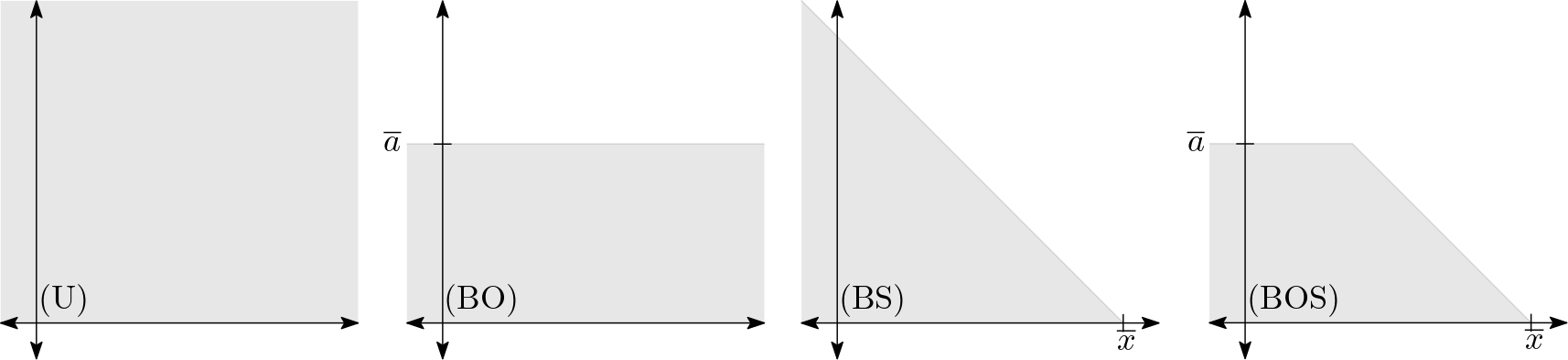}
	\caption{The graph of the feasible actions multifunction $A$ under the four inventory models (U, BO, BS, BOS) with backorders, respectively.}
	\label{fig:inventory_models}
\end{figure}

The inventory model can be formulated as an MDP in the following way. Let the state space $\XX$ represent the amount of available inventory and the action space $\AA := [0,+\infty)$ represent order sizes. The multifunction of feasible actions $A : \XX \to 2^{\AA} \setminus \{\emptyset\}$, where $A(x) \subset \AA$ represents the possible inventory order sizes at the state $x \in \XX$ and captures the constraints implied by the values of $\bar a$ and $\bar x$. With each combination of the pair $(\bar a, \bar x)$, there are four models characterized by the state space $\XX$ and feasible actions $A(x)$, listed below.
\begin{enumerate}
    \item (U) Unbounded orders, unbounded storage: $\XX := \RR$ for problems with backorders and $\XX := [0,+\infty)$ for problems with lost sales, and $A(x) := [0,+\infty)$ for each $x \in \XX$.
    \item (BO) Bounded orders, unbounded storage: $\XX := \RR$ for problems with backorders and $\XX := [0,+\infty)$ for problems with lost sales, and $A(x) := [0,\bar a]$ for each $x \in \XX$.
    \item (BS) Unbounded orders, bounded storage: $\XX := (-\infty, \bar x]$ for problems with backorders and $\XX := [0,\bar x]$ for problems with lost sales, and $A(x) := [0, 0 \lor (\bar x - x)]$ for each $x \in \XX$.
    \item (BOS) Bounded orders, unbounded storage: $\XX := (-\infty, \bar x]$ for problems with backorders and $\XX := [0,\bar x]$ for problems with lost sales, and $A(x) := [0,\bar a] \cap [0, 0 \lor (\bar x - x)]$ for each $x \in \XX$.
\end{enumerate}
In the unbounded orders, unbounded storage model (U), $A$ is not compact valued. In the bounded storage, unbounded orders model (BO), the image $A(\XX) = \bigcup_{x \in \XX} A(x) = (-\infty, \bar x]$ is not a compact set in the model with backorders, and $A(\XX) = [0, \bar x]$ in the model with lost sales. In the bounded orders, unbounded storage model (BS) and the bounded orders, bounded storage model (BOS), $A$ is compact valued and has a compact image. In each model, $A$ is upper semicontinuous and lower semicontinuous as a multifunction, and the graph of $A$ is closed. The graphs of the feasible actions $A$ in each of the models (U, BO, BS, BOS) for models with backorders are given in Figure~\ref{fig:inventory_models}.

The dynamics of the system are defined by the equation
\begin{align}
	x_{t+1} &= T(x_t + a_t - D_{t+1}), && t = 0, 1, \dots,
  \label{eq:transition_dynamics}
\end{align}
where $x_t$ and $a_t$ denote the current inventory level and the amount ordered at period $t$. The transition probability $q(dx_{t+1} | x_t, a_t)$ for the MDP defined by the dynamics above is
\begin{equation}
	q(B | x_t, a_t) = P(T(x_t + a_t - D_{t+1}) \in B),
\end{equation}
for each measurable $B \subset \RR$. The one-step expected cost is
\begin{align}
	c(x,a) &:= K \one{a > 0} + \bar c a + \EE h(T(x+a-D)), && x \in \XX, a \in \AA,
\end{align}
where $\one*{B}$ is the indicator of the set $B$. We observe that $c$ is inf-compact.

Let $H_t = (\XX \times \AA)^t \times \XX$ be the space of histories of the system up to period $t = 0, 1, \dots$. A (possibly randomized) control at $t = 0, 1, \dots$ is a transition probability $\phi_t : H_t \to \PP(\AA)$. A policy $\phi$ is a sequence $(\phi_0, \phi_1, \dots)$ of controls. The policy $\phi$ is Markov if all controls $\phi_t$ depend only on the current state and time. A Markov policy is stationary if all controls only depend on the state (and not time). A policy $\phi$ is nonrandomized if each control $\phi_t$ is concentrated at a single point in $\AA$. In this latter case, there is no loss in generality to write $\phi_t : \XX \to \AA$.

For an horizon $N = 1, 2, \dots, +\infty$, the expected total discounted costs is defined as
\begin{align}
	v^\pi_{N, \alpha}(x) &:= \EE_x^\pi \left[ \sum_{t=0}^{N-1} \alpha^t c(x_t, a_t) \right], && x \in \XX,
	\label{eq:discounted_value_fn}
\end{align}
where $v^\pi_{0, \alpha}(x) := 0$ ($x \in \XX$). When $N=+\infty$,~\eqref{eq:discounted_value_fn} defines the infinite-horizon expected total discounted cost, which is denoted by $v^\pi_\alpha(x)$. Let $v_{N,\alpha} := \inf_{\pi} v_{N,\alpha}^\pi(x)$ and $v_\alpha := \inf_{\pi} v_\alpha^\pi(x)$, $x \in \XX$. A policy $\pi$ is called \emph{$N$-stage optimal} if $v^\pi_{N,\alpha} = v_{N,\alpha}$ and \emph{discount-optimal} if $v^\pi_\alpha = v_\alpha$. The optimality equations
\begin{align}
	v_{t+1, \alpha}(x)
	&= \min \set{G_{t,\alpha}(x), K + \min_{a \in A(x)} G_{t,\alpha}(x+a)} - \bar c x
	&& x \in \XX,\quad t=0, 1, \dots,
	\label{eq:finite_horizon_bellman} \\
	v_{\alpha}(x)
	&= \min \set{G_{\alpha}(x), K + \min_{a \in A(x)} G_{\alpha}(x+a)} - \bar c x
	&& x \in \XX,
	\label{eq:infinite_horizon_bellman}
\end{align}
where
\begin{align}
	G_{t,\alpha}(x)
  &:= \bar c x + \EE h(T(x-D)) + \alpha \EE v_{t,\alpha}(T(x-D)),
  && x \in \XX, \quad t = 0, 1, \dots,
	\label{eq:finite_horizon_g} \\
	G_{\alpha}(x)
  &:= \bar c x + \EE h(T(x-D)) + \alpha \EE v_{\alpha}(T(x-D)),
  && x \in \XX,
	\label{eq:infinite_horizon_g}
\end{align}
hold in each model (U, BO, BS, BOS); see Section~\ref{sec:continuity} for details.

The following theorem shows that both $v_{N,\alpha}$ and $G_{N,\alpha}$ are continuous for each $N=0,1,\dots$. For the model (U) with backordering,~\citet[Theorem 8.3.4]{simchi_levi_2014_logic} proved continuity when the holding costs are linear, and~\citet[Theorem 5.2]{feinberg_2017_structure} proved continuity when the holding costs are convex. When the models (U, BO, BS, BOS) are considered together with backorders and with lost sales, continuity of values is a corollary of a generalized form of Berge's maximum theorem, which holds for discontinuous one-step cost functions. \citet[Theorem 14]{feinberg_2021_continuity} states the continuity of the functions $v_{N, \alpha}$ for the models (U, BO, BS, BOS), and its proof established continuity of  $x \mapsto \EE v_{N, \alpha}(x-D);$ see additional references there. The next theorem restates~\citet[Theorem 14]{feinberg_2021_continuity} with the additional claim that $G_{N,\alpha}$ are also continuous.

\begin{theorem}[{cf.~\citet[Theorem 14]{feinberg_2021_continuity}}]
	For each inventory model (U, BO, BS, BOS) with lost sales and with backorders, and for each $N=0,1,\dots$, the functions $v_{N,\alpha}$, $x \mapsto \EE v_{N, \alpha}(T(x-D))$, and $G_{N, \alpha}$ are continuous, and there is an optimal $N$-stage deterministic Markov policy $\phi^N = (\phi_0, \dots, \phi_{N-1})$.
	\label{thm:continuity_of_finite_horizon_values}
\end{theorem}
The following theorem shows that the discounted value function $v_\alpha$ is continuous. For the model (U) with backordering,~\citet[Theorem 9.11]{bensoussan_2011_dynamic} proved continuity when the holding costs are linear, and~\citet[Theorem 5.3]{feinberg_2017_structure} proved continuity when the holding costs are convex. Continuity under the models (BO, BOS) with lost sales and with backorders is a new result. Furthermore, the following theorem establishes continuity without constructing discount-optimal policies. Instead, the result follows from the continuity of discounted value function in the problem with $K=0$.

\begin{theorem}
	For each inventory model (U, BO, BS, BOS) with lost sales and with backorders, the functions $v_\alpha$, $x \mapsto \EE v_\alpha(T(x-D))$, and $G_\alpha$ are continuous. For the models (U, BS), each of these functions is finite. For the inventory models (BO, BOS), if any of these functions is finite at a point $x_0 \in \XX$, then they are each finite for all $x \in \XX$.

	\label{thm:continuity_of_infinite_horizon_values}
\end{theorem}

In Section~\ref{sec:continuity} we provide an example for the models (BO, BS) where the discounted value function is uniformly infinite, i.e., $v_\alpha(x) = +\infty$ for each $x \in \XX$; see Example~\ref{ex:infinite_values} for details.
\revision{%
  The remainder of this section concerns the structure of optimal policies for the inventory models (U, BS). The next definition is standard in the inventory control literature; see, e.g.,~\citet[Definition 9.2.1]{simchi_levi_2014_logic}.
}
We recall that a function $f : \RR \to \RR$ is \emph{$K$-convex} on $X$ if for all $x, y \in X$ with $x \leq y$, and for all $\theta \in [0,1]$, the inequality
\begin{equation}
  f(\theta x + (1-\theta)y) \leq \theta f(x) + (1-\theta) f(y) + (1-\theta) K.
\end{equation}
holds. For example, every convex function is $0$-convex. For inf-compact $K$-convex functions, the quantities
\begin{align}
  S &\in \argmin \set{ f(x) : x \in \XX}, \label{eq:big_S} \\
  s &= \min \set{x \in \XX : f(x) \leq K + f(S)}, \label{eq:little_s}
\end{align}
exist, where $s \leq S$ by construction. The numbers $s$ and $S$ are used to define optimal $(s,S)$ policies.

\begin{definition}[$(s,S)$ policies]
  For each $t = 0, 1, \dots$, let $s_t \leq S_t$ be finite numbers. An $(s_t, S_t)$ control at step $t$ is the function $\phi_t(x_t) = (S_t - x_t) \one{x_t < s_t}$. A Markov policy $\phi$ is an $(s_t, S_t)$ policy if $\phi = (\phi_0, \phi_1, \dots)$ are all $(s_t,S_t)$ controls for each $t$. A policy is called an $(s, S)$ policy if it is a stationary Markov $(s_t, S_t)$ policy.
\end{definition}
Determination of the values of $s$ and $S$ follows the classic analysis of $K$-convex functions.
\revision{%
For some finite-horizon problems, an $(s_{t}, S_{t})$ policy is optimal. For other problems, the policy which never orders is optimal. However, there are problems where it is optimal to play an $(s_{t}, S_{t})$ policy up to a certain step, but never to order afterward. Abridged $(s_t,S_t)$ policies of this type are given in the next definition.
\begin{definition}[$(s_t,S_t,n)$ policies]
  A finite horizon policy $\phi = (\phi_0, \phi_1, \dots, \phi_{N-1})$ is an $(s_t, S_t, n)$ policy if $\phi_t$ is an $(s_t, S_t)$ control for each $t = 0, 1, \dots, (N-n-1)$, and $\phi_t(x) \equiv 0$ for each $t = N-n, \dots, N-1$.
\end{definition}
We observe that an $(s_{t}, S_{t}, 0)$ policy corresponds to a classic $(s_{t}, S_{t})$ policy. Therefore $(s_t, S_t, n)$ policies generalize finite-horizon $(s_t, S_t)$ policies.
}%
\revision{%
The discount-optimality of stationary $(s,S)$ policies holds for sufficiently large discount factors $\alpha < 1$. We denote total demand over $t$ periods by $\mathbf{S}_{t} := \sum_{i=1}^{t} D_i$. Following~\citet[(4.1)-(4.3)]{feinberg_2017_structure}, denote
\begin{align}
    \alpha^\star &:= 1 + \lim_{x \to -\infty} \frac{h(x)}{\bar cx}, \label{eq:alpha_star}
\end{align}
where $-\infty \leq \alpha^\star < 1$. For each $t = 0, 1, \dots$ and $\alpha \in [0, 1)$ define the function
\begin{align}
  f_{t,\alpha}(x) := \bar c x + \sum_{i=0}^{t} \alpha^i \EE h(x - \mathbf{S}_{i+1}), \qquad x \in \XX.
  \label{eq:f_t_alpha}
\end{align}
We observe that $f_{t,\alpha}$ is convex for each $t$ and $\alpha$. Consider the number
\begin{align}
  N_\alpha &:= \inf \set*{t \in \NN : \lim_{x \to -\infty} f_{t,\alpha}(x) = +\infty}, \label{eq:N_alpha}
\end{align}
where $\inf \emptyset := +\infty$. Since $h$ is nonnegative, the number $N_\alpha$ is a nonincreasing function of $\alpha$.~\citet[Theorem 4.3]{feinberg_2017_structure} describe the structure of optimal policies for the model (U) with backorders for all discount factors $\alpha \in [0,1)$. The following theorem extends this result to the model (BS) with backorders.
\begin{theorem}
  For the inventory models (U, BS) with backorders, consider the discount factor $\alpha \in [0,1)$.
\begin{enumerate}
 	\item If $\alpha \leq \alpha^\star$, then the following statements hold.
	\begin{enumerate}
	    \item For each $N = 1, 2, \dots$, the policy that never orders is $N$-stage optimal.
	    \label{item:discounted_finite_no_order}
	    \item The policy that never orders is discount optimal. \label{item:discounted_infinite_no_order}
	\end{enumerate}
	  \item If $\alpha^\star < \alpha$, then the following statements hold.
	\begin{enumerate}
	    \item The equality
	    \begin{equation}
	        N_\alpha = \min \set*{t \in \NN :
	            \frac{\alpha^\star}{1-\alpha^\star} < \sum_{i=1}^{t} \alpha^i
	        }
        \label{eq:N_alpha_formula}
	    \end{equation}
	    holds, and $0 \leq N_\alpha < +\infty$. \label{item:N_alpha_equation}
	    \item For each $N = 1, 2, \dots$,
	    \begin{enumerate}
	        \item if $N \leq N_\alpha$, then the policy that never orders is $N$-stage optimal;
	        \item if $N > N_\alpha$, then the policy $(s_{t,\alpha}, S_{t,\alpha}, N_\alpha)$, where $S_{t,\alpha} := S_{G_{N-t-1, \alpha}}$ and $s_{t,\alpha} := s_{G_{N-t-1,\alpha}}$ are defined in~\eqref{eq:big_S} and~\eqref{eq:little_s} for $t=0,\dots, N-1$, is $N$-stage optimal. \label{item:alpha_star_lt_alpha_N_gt_N_alpha}
	    \end{enumerate}
	    \label{item:discounted_finite_s_S}
	    \item For the pair $(s_\alpha, S_\alpha)$ defined in \eqref{eq:big_S} and \eqref{eq:little_s} for the function $f(x) := G_\alpha(x)$, the $(s_\alpha, S_\alpha)$ policy is discount optimal. In addition, the sequence $\set{(s_{N,\alpha}, S_{N,\alpha})}_{N=N_\alpha+1, N_\alpha+2, \dots}$ defined in Statement~\ref{item:alpha_star_lt_alpha_N_gt_N_alpha} is bounded, and each of its limit points $(s_\alpha', S_\alpha')$ is discount optimal. \label{item:discounted_infinite_s_S}
	\end{enumerate}
\end{enumerate}

\label{thm:bounded_storage_s_S_policies}
\end{theorem}
}

\revision{%
  The case $\alpha^\star < 0$ deserves special attention. According to~\citet[Lemma 4.1]{feinberg_2017_structure}, this condition is equivalent to the existence of inventory levels $z, y \in \XX$ ($z < y$) satisfying
  \begin{equation}
      \frac{h(y)-h(z)}{y-z} < -\bar c,
  \end{equation}
  which for the model (U) with backorders is a well-known condition in the literature; see, e.g.,~\citet{scarf_1960_optimality,iglehart_1963_optimality}, and~\citet{veinott_1965_computing}. The following corollary summarizes Theorem~\ref{thm:bounded_storage_s_S_policies} in the special case of $\alpha^\star < 0$.

  \begin{corollary}
    For the inventory models (U, BS) with backorders and discount factor $\alpha \in [0,1)$, suppose that $\alpha^\star < 0$ holds. Then $N_\alpha = 0$, and the following statements hold.
      \begin{enumerate}
          \item For each $N = 1, 2, \dots$, the policy $(s_{t,\alpha}, S_{t,\alpha})$, where  $S_{t,\alpha} := S_{G_{N-t-1, \alpha}}$ and $s_{t,\alpha} := s_{G_{N-t-1,\alpha}}$ are defined in~\eqref{eq:big_S} and~\eqref{eq:little_s} for $t=0,\dots, N-1$, is $N$-stage optimal.
          \label{item:alpha_star_negative_finite_horizon}
          \item For the pair $(s_\alpha, S_\alpha)$ defined in \eqref{eq:big_S} and \eqref{eq:little_s} for the function $f(x) := G_\alpha(x)$, the $(s_\alpha, S_\alpha)$ policy is discount optimal. In addition, the sequence $\set{(s_{N,\alpha}, S_{N,\alpha})}_{N=1,2, \dots}$ defined in Statement~\ref{item:alpha_star_negative_finite_horizon} is bounded, and each of its limit points $(s_\alpha', S_\alpha')$ is discount optimal. \label{item:alpha_star_negative_infinite_horizon}
      \end{enumerate}
  \end{corollary}

  We close this section with the following remark about average-cost optimal $(s,S)$ policies.~\citet{schal_1993_average} considers the model (BS) with backorders, inf-compact holding costs $h(x)$, and left continuous, nondecreasing order costs. In Theorem 4.8 there,~\citeauthor{schal_1993_average} showed that (i) if there is a sequence of discount factors $\alpha_k \uparrow 1$ such that an $(s_{\alpha_k}, S_{\alpha_k})$ policy is $\alpha_k$-discount optimal, then there exists a stationary average-cost optimal policy; (ii) if there exists $\alpha_0 < 1$ such that $(s_\alpha, S_\alpha)$ policies are $\alpha$-discount optimal for all $\alpha \in [\alpha_0, 1)$, then there exists an average-cost optimal $(s,S)$ policy. In our Theorem~\ref{thm:bounded_storage_s_S_policies}, we showed that for the model (BS) with convex holding costs and order costs of the form $a \mapsto (K \one{a > 0} + \bar c a)$, there exist $\alpha$-discount optimal $(s_\alpha, S_\alpha)$ policies for all $\alpha \in [\alpha^\star, 1)$. Since $\alpha^\star < 1$,~\citet[Theorem 4.8]{schal_1993_average} implies that there exists an average-cost optimal $(s,S)$ policy for the model (BS).~\citet[Theorem 4.4]{feinberg_2017_optimality} proves the average-cost optimality of $(s,S)$ policies for the model (U) with backorders and convex holding costs. For the model (U) the similar result is proved in~\citet{feinberg_2017_structure} by using the results on the existence of stationary optimal policies for average-cost MDPs with noncompact action sets~\cite{feinberg_2012_average}.
}

\begin{figure}[t]
  \centering
  \includegraphics[scale=0.25]{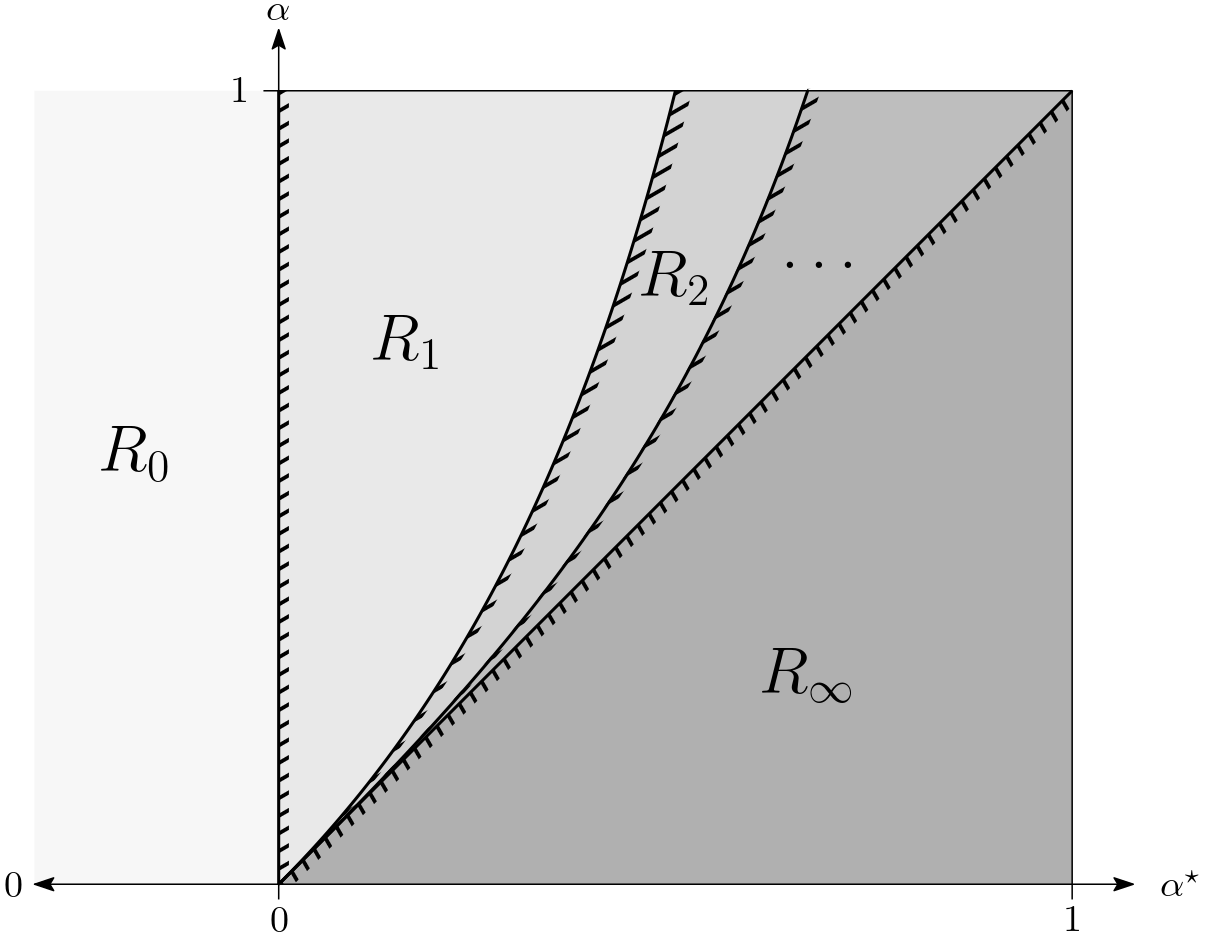}
  \begin{center}
    \begin{tabular}{@{}cclll@{}} \toprule
      \multicolumn{4}{@{}l@{}}{Region} & \multicolumn{1}{l@{}}{Structure of optimal policy} \\ \midrule
      $R_0$ & $=$ & $[-\infty, 0) \times [0,1)$ && $(s_t,S_t)$\\
      $R_1$ & $=$ & $\set{(\alpha^\star, \alpha) \in [0,1)^2 : \frac{\alpha^\star}{1-\alpha^\star} < \alpha}$ &&  $(s_t,S_t,1)$ \\
      $R_n$ & $=$ & $\set{(\alpha^\star, \alpha)\in [0,1)^2 :\frac{\alpha^\star}{1-\alpha^\star} < \sum_{i=1}^{n} \alpha^i} \setminus \bigcup_{i=1}^{n-1} R_{i}$ & $n=2, 3, \dots$ & $(s_t,S_t,n)$\\
      $R_{\infty}$ &$=$& $\set{(\alpha^\star, \alpha) \in [0,1)^2 : \alpha \leq \alpha^\star}$ && Never-order \\\bottomrule
    \end{tabular}
  \end{center}
  \caption{%
    Structure of optimal policies for a discounted $N$-horizon problem with $\alpha \in [0,1)$. In $R_0$, an $(s_{t}, S_{t})$ policy is optimal. In $R_n$ ($n = 1, 2, \dots$), an $(s_{t}, S_{t}, n)$ policy is optimal. In $R_{\infty}$, the policy which never orders is optimal.
  }
  \label{fig:structure_finite_horizon}
\end{figure}

\begin{figure}[t]
  \centering
   \includegraphics[scale=0.25]{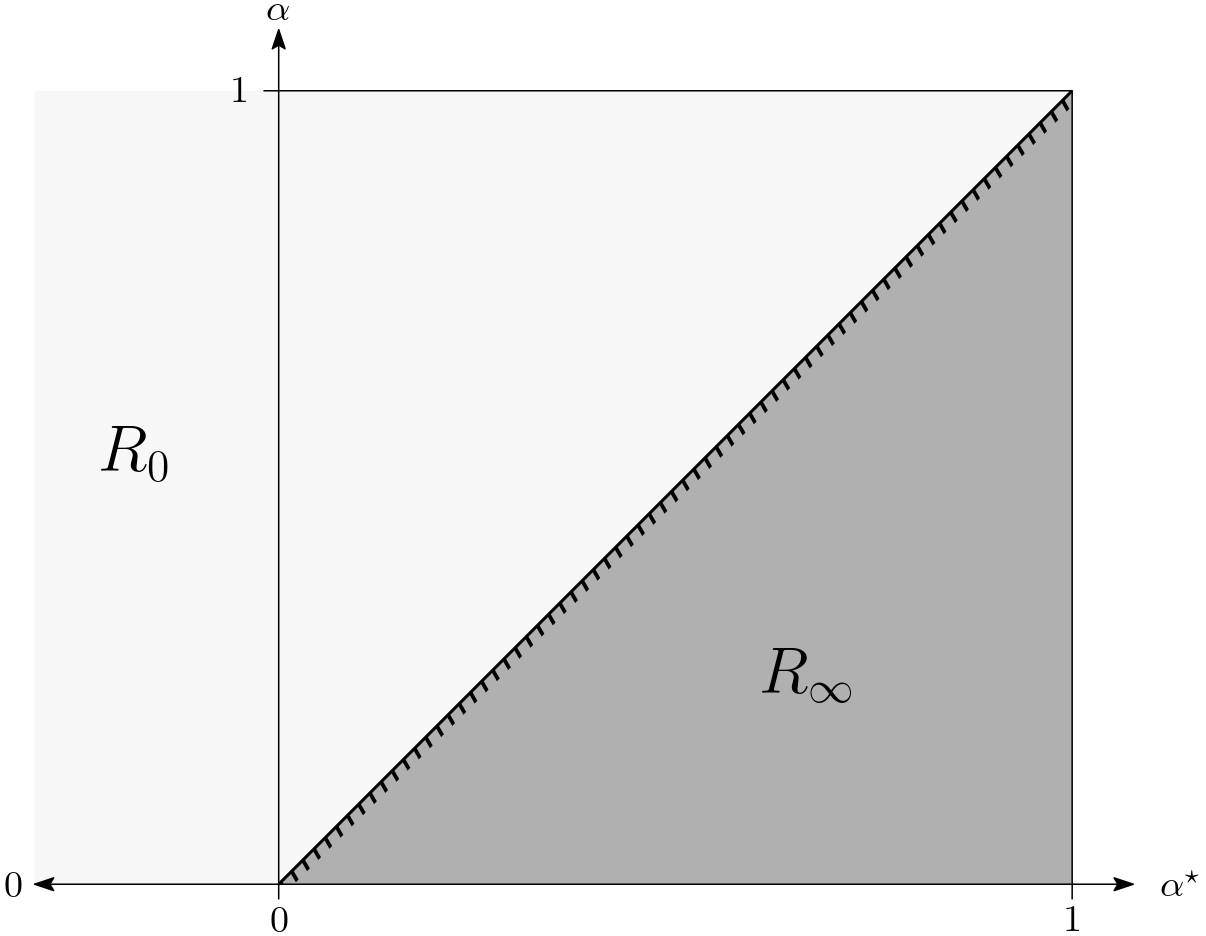}
  \begin{center}
    \begin{tabular}{@{}ccll@{}l@{}} \toprule
      \multicolumn{4}{@{}l@{}}{Region} & \multicolumn{1}{@{}l@{}}{Structure of optimal policy} \\ \midrule
      $R_0$ & $=$ & $\set{(\alpha^\star, \alpha) \in [-\infty, 1) \times [0, 1): \alpha^\star < \alpha}$ && $(s,S)$\\
      $R_{\infty}$ &$=$& $\set{(\alpha^\star, \alpha) \in [0,1)^2: \alpha \leq \alpha^\star}$ && Never-order \\\bottomrule
    \end{tabular}
  \end{center}
  \caption{Structure of optimal policies for a discounted infinite-horizon problem with $\alpha \in [0,1)$. In $R_0$, an $(s, S)$ policy is optimal. In $R_\infty$, the policy which never orders is optimal.}
  \label{fig:structure_infinite_horizon}
\end{figure}
\section{%
  Proofs of Theorems~\ref{thm:continuity_of_finite_horizon_values},~\ref{thm:continuity_of_infinite_horizon_values}, and Example
}
\label{sec:continuity}

\revision{%
In this section we establish the correctness of Theorems~\ref{thm:continuity_of_finite_horizon_values} and~\ref{thm:continuity_of_infinite_horizon_values}.
}
For a function $f : X \to \RR$, we denote \revision{the domain of $f$ by
\begin{equation}
	\dom f = \set{x \in X: f(x) < +\infty}
\end{equation}
and}, for each $\lambda \in \RR$ \revision{and subset $C \subset X$}, the sublevel set of $f$ restricted to $C$ by
\begin{equation}
	\lev_f(\lambda; C) = \set{x \in C : f(x) \leq \lambda},
\end{equation}
where $\lev_f (\lambda) := \lev_f(\lambda; X)$. 
For a multifunction $\Gamma : X \to Y$ and subset $K \subset X$, the graph of $\Gamma$ on $K$ is denoted by $\Gr_K(\Gamma) := \set{(x,y) : x \in K , y \in \Gamma(x)}$, and $\Gr(\Gamma) := \Gr_X(\Gamma)$.
The next definition was introduced in~\citet{feinberg_2012_average} in the conditions of Assumption~\ref{ass:w_star}, but was formulated and studied as a separate property in the later works~\citet{feinberg_2013_berge, feinberg_2014_berge, feinberg_2015_continuity,feinberg_2021_continuity} with applications relevant to inventory control in~\citet{feinberg_2016_optimality}.
\begin{definition}[\revision{\citet[Definition 1.1]{feinberg_2013_berge}}]
	For metric spaces $X$ and $Y$, a multifunction $\Gamma : X \to 2^Y \setminus \
	\{\emptyset\}$, and function $f : X \times Y \to \RR \cup \set{+\infty}$, the function $f$ is \emph{$\KK$-inf-compact} on $\Gr(\Gamma)$ if, for each nonempty compact set $K \subset X$ and $\lambda \in \RR$, the set $\lev_f(\lambda; K \times Y)$ is compact.
\end{definition}

The following assumption, known as Assumption~(\ref{ass:w_star}), implies the correctness of the optimality equations for finite and infinite-horizon MDPs.

\begin{specialAssumption}[\setword{W*}{ass:w_star}]
	The following hold.
	\begin{enumerate}
		\item The function $c$ is $\KK$-inf-compact on $\Gr(A)$ and bounded below;
		\item the transition probability $q(\:\cdot\:|x,a)$ is weakly continuous in $(x,a) \in \Gr(A)$.
	\end{enumerate}
\end{specialAssumption}

It is well known that the model (U) with backorders satisfies Assumption~(\ref{ass:w_star}); see, e.g.~\cite{heyman_2004_stochastic, bertseks_2005_dynamic, simchi_levi_2014_logic, feinberg_2016_optimality}. In fact, Assumption~(\ref{ass:w_star}) is satisfied by each model (U, BO, BS, BOS) with backorders and with lost sales. In each model, $c$ is $\KK$-inf-compact because the set $\Gr(A)$ is closed, which implies that
\begin{equation}
	\lev_f(\lambda; \Gr_K(A))
	= \set{(x,a) : x \in K, a \in A(x), c(x,a) \leq \lambda}
	\subset K \times [0, \lambda / \bar c]
\end{equation}
for each compact $K \subset \XX$ and $\lambda \in \RR$. Further, in each model the transition dynamics~\eqref{eq:transition_dynamics} is a continuous function of the state, action, and demand; so it follows (see, e.g., \citet[p. 92]{hernandez-lerma_1989_adaptive}) that $q(\:\cdot\:|x,a)$ is weakly continuous. Therefore, the optimality equations~\eqref{eq:finite_horizon_bellman} and~\eqref{eq:infinite_horizon_bellman} are justified according to \citet[Theorem 2]{feinberg_2012_average},
which also shows that each $v_{t,\alpha}$ is lower semicontinuous and $v_{t,\alpha} \uparrow v_\alpha$. Furthermore, since $v_{t,\alpha} \uparrow  v_{\alpha}$ implies $G_{t,\alpha} \uparrow G_{\alpha}$, and since $v_{t,\alpha}(x) \geq \EE h(x-D) \geq h(x - \EE D)$ for each $x \in \XX$ and $t \geq 1$,  it follows that
\begin{align}
	\lim_{x \to +\infty} G_{t, \alpha}(x) = \lim_{x \to +\infty} v_{t, \alpha}(x) =
	\lim_{x \to +\infty} G_{\alpha}(x) = \lim_{x \to +\infty} v_{\alpha}(x) = +\infty,
	\label{eq:g_right_limit}
\end{align}
because $\lim_{\abs{x} \to +\infty} h(x-\EE D) = +\infty$.
\begin{proof}[Proof of Theorem~\ref{thm:continuity_of_finite_horizon_values}]
  Let $N \in \NN$ be arbitrary. The continuity of $v_{N,\alpha}$ and existence of $\phi^N$ follow directly from the statement of Theorem 14 in~\citet{feinberg_2021_continuity}. The proof of that theorem also establishes that $x \mapsto \EE v_{N,\alpha}(T(x-D))$ is continuous. Then the function $G_{N, \alpha}$ defined in~\eqref{eq:finite_horizon_g} is the sum of continuous functions; hence, continuous.
\end{proof}

\revision{
  The next objective is to show that $v_\alpha$ is continuous. We observe that for the setup-cost inventory control problem with $K = 0$, the one-step expected cost is convex. Since $0 \leq K \one{a > 0}$, the controller will always incur a one-step cost \emph{no less than} they would for the same order in the problem without setup costs. But since $K \one{a > 0} \leq K$, they will always incur a one-step cost \emph{no greater than} $K$ plus that same no-setup cost amount. The upshot is that the value function $v_\alpha$ can be estimated above and below by two convex functions. This discussion is formalized in Lemma~\ref{lem:no_setup_costs}, which is exploited by Theorem~\ref{thm:continuity_of_infinite_horizon_values} to show that $v_\alpha$ is continuous.

  The following technical lemma provides sufficient conditions for the parametric minimum of a convex function of two variables to be convex, which is needed for Lemma~\ref{lem:no_setup_costs}. The lemma follows from two facts: Berge's theorem (see, e.g.,~\citet[Theorem 6]{feinberg_2021_continuity}), which implies lower semicontinuity of the value function and existence of solutions; and from~\citet[Proposition 2.22(a)]{rockafellar_2009_variational}, which implies convexity of the value function and solution sets in an unconstrained problem when the objective function is convex. The next lemma extends~\cite[Proposition 2.22(a)]{rockafellar_2009_variational} to constrained parametric optimization problems.

  \begin{lemma}
    Let $X$ and $Y$ be Euclidean spaces, let $\Phi : X \to 2^{Y}$ be a multifunction such that $\Gr(\Phi) \subset X \times Y$ is convex, and let $u : X \times Y \to \RR$ be $\KK$-inf-compact on $\Gr(\Phi)$ and convex. Then
    \begin{align}
      u^*(x) &:= \inf_{y \in \Phi(x)} u(x,y), &&x \in X,
    \end{align}
    is convex and lower semicontinuous, and $\Phi^*(x) := \argmin_{y \in \Phi(x)} u(x,y)$ is nonempty, upper semicontinuous, and compact and convex-valued.
    \label{lem:bmt_convex_functions}
  \end{lemma}

  \begin{proof}
    That $u^*$ is lower semicontinuous and that $\Phi^*$ is upper semicontinuous with compact values follows from~\citet[Theorem 5, see also Theorem 13]{feinberg_2021_continuity}. It is sufficient, then, to show that $u^*$ is convex and that $\Phi^*$ has convex values. In fact, the result is an immediate corollary of~\citet[Proposition 2.22(a)]{rockafellar_2009_variational}, after we perform the following modification. Let $\psi$ denote the (convex) indicator function of the set $\Gr(\Phi) \subset X \times Y$. That is, $\psi(x,y) = 0$ if $y \in \Phi(x)$ and $\psi(x,y) = +\infty$ otherwise. Since $\Gr(\Phi)$ is convex, the function $\psi$ is convex.
    Consider the unconstrained problem $\tilde{u}(x,y) := u(x,y) + \psi(x,y)$ with $\tilde{\Phi}(x) := Y$ for each $x \in X$.
    Since the equality
    \begin{equation*}
      u^*(x) = \inf_{y \in \Phi(x)} u(x,y) = \inf_{y \in Y} \set{u(x,y) + \psi(x,y)} = \tilde{u}^*(x)
    \end{equation*}
    holds, it is sufficient to determine if $\tilde{u}$ is convex and $\tilde{\Phi}^* = \Phi^*$ has convex values. This follows from~\citet[Proposition 2.22(a)]{rockafellar_2009_variational}.
  \end{proof}

  The next lemma establishes the relationship between the discounted values $v_\alpha$ for the inventory problem with setup costs and the discounted values $v_\alpha^0$ for the corresponding problem without setup costs.

  \begin{lemma}
    For each inventory model (U, BO, BS, BOS) with lost sales and with backorders, consider the problem without setup costs, i.e., with one-step costs
    \begin{align}
      c^0(x,a) &:= \bar c a + \EE h(T(x+a-D)), && x \in \XX, a \in \AA,
      \label{eq:no_setup_cost}
    \end{align}
    with finite horizon values $v_{N,\alpha}^0$ and discounted values $v_{\alpha}^0$. The following statements hold.
    \begin{enumerate}
      \item For each $N = 0, 1, \dots$, the function $v_{N, \alpha}^0$ is finite, convex, continuous, and
      \begin{align}
        v_{N,\alpha}^0(x) \leq v_{N, \alpha}(x) \leq v_{N, \alpha}^0(x) + \frac{1-\alpha^N}{1-\alpha}K, && x \in \XX.
        \label{eq:no_setup_inequality_finite_horizon}
      \end{align}
      \item The function $v_{\alpha}^0$ is convex, continuous, and
      \begin{align}
        v_{\alpha}^0(x) \leq v_{\alpha}(x) \leq v_{\alpha}^0(x) + \frac{K}{1-\alpha}, && x \in \XX.
        \label{eq:no_setup_inequality_infinite_horizon}
      \end{align}
      Furthermore, if $v_\alpha^0(x_0) < +\infty$ for some $x_0 \in \XX$, then $v_{\alpha}^0(x) < +\infty$ for all $x \in \XX$.
    \end{enumerate}
    \label{lem:no_setup_costs}
  \end{lemma}

  \begin{proof}
    The function $v_{0,\alpha}^0 \equiv 0$ is convex by definition. Suppose $v_{N, \alpha}^0$ is convex, finite, and continuous.
    Since the function $c^0$ is jointly convex in $x$ and $a$, Lemma~\ref{lem:bmt_convex_functions} implies that
    \begin{equation*}
      v_{N+1,\alpha}^{0}(x) = \min_{a \in A(x)} \set{c^0(x,a) + \alpha \EE v_{N,\alpha}^0(T(x+a-D))}
    \end{equation*}
    is a convex function. Since from Theorem~\ref{thm:continuity_of_finite_horizon_values} each $v_{N,\alpha}$ is finite, the inequalities
    \begin{align*}
      v_{N+1,\alpha}^0(x) &= \min_{a \in A(x)} \set{c^0(x,a) + \alpha \EE v_{N,\alpha}^0(T(x+a-D))} \\
      &\leq \min_{a \in A(x)} \set{c(x,a) + \alpha \EE v_{N,\alpha}(T(x+a-D))} \\
      &= v_{N+1,\alpha}(x) \\
      &\leq \min_{a \in A(x)} \set{K + c^0(x,a) + \alpha \EE \left[v_{N,\alpha}^0(T(x+a-D)) + (1-\alpha^t)(1-\alpha)^{-1}\right]} \\
      &= v_{N+1,\alpha}^0(x) +  (1-\alpha^{N+1})(1-\alpha)^{-1}K.
    \end{align*}                                                                imply that $v_{N+1, \alpha}^0(x)$ is finite for each $x \in \XX$; hence, also continuous. Since $v_{N,\alpha} \uparrow v_{\alpha}$ and $v_{N,\alpha}^0 \uparrow v_{\alpha}^0$, it follows that
    \begin{equation*}
      v_{\alpha}^0 \leq v_{\alpha} \leq v_{\alpha}^0 + (1-\alpha)^{-1}K.
    \end{equation*}

    We now show that $v_\alpha^0$ is continuous. We observe that $v_\alpha^0$ is convex, lower semicontinuous, and nonnegative. If $v_\alpha^0 \equiv +\infty$, continuity is evident. Instead, suppose there exists $x_0 \in \dom v_\alpha^0$. We shall show that $\dom v_\alpha^0 = \XX$, which will imply that $v_\alpha^0$ is continuous.

    Suppose $\dom v_\alpha^0 = \set{\tilde{x}}$. Then there exists $a \in A(\tilde{x})$ such that $a = D$ almost surely. The multifunction $A$ has the property that $A(x) \supset A(y)$ if $x \leq y$. Therefore, $a \in A(x)$ for all $x \in \XX$ satisfying $x \leq \tilde{x}$; whence $v_\alpha^0(x) \leq (\bar c a + \EE h(T(x+a -D)))(1-\alpha)^{-1} < +\infty$ for all such $x$, a contradiction. Therefore, $\dom v_\alpha^0$ is an interval of positive length.

    To show that $\dom v_\alpha^0 = \XX$, consider the left segment $\XX_\ell := \XX \cap (-\infty, \tilde{x}]$ and the right segment $\XX_r := \XX \cap [\tilde{x}, +\infty)$. We show that $\XX_\ell \subset \dom v_\alpha^0$. There are two cases.
    \begin{enumcase}
        \item $\XX_\ell = (-\infty, \tilde{x}]$. This is the model with backorders. Let $x_\ell = \inf \dom v_\alpha^0$. Suppose for the sake of contradiction that $x_\ell > -\infty$. Let $\phi^*$ denote an optimal policy. Then for all $\varepsilon > 0$ we can fix $x_\varepsilon \in \dom v_\alpha^0$ such that $x_\varepsilon - x_\ell < \varepsilon$ and such that $P(T(\varepsilon + \phi^*(x_\varepsilon) > D))$. We consider two cases. If $\varepsilon + \phi^*(x_\varepsilon) \in A(x_\ell)$, then we let $a_\ell := \varepsilon + \phi^*(x_\varepsilon)$. If, however, $\varepsilon + \phi^*(x_\varepsilon) \notin A(x_\ell)$, then we let $a_\ell := \limsup_{\varepsilon \to 0} \phi^*(x_\varepsilon)$, and we observe that $a_\ell \in A(x_\ell)$, since $A(x_\varepsilon) \subset A(x_\ell) \subset [0, a_\varepsilon]$ for each $\varepsilon > 0$. In either case, $a_\ell$ has the property that, under the policy $\phi(x) := a_\ell \one{x \leq x_\ell}$ starting at $x \leq x_\ell$, the system never leaves the interval $[x, x + a_\ell]$. Then $v_\alpha^0(x) \leq (\bar c a_\ell + 2\EE h(T(x+a_\ell-D)))(1-\alpha)^{-1} < +\infty$, a contradiction. Therefore, $(-\infty, \tilde{x}] \subset \dom v_\alpha^0$.
        \item $\XX_\ell = [0, \tilde{x}]$. This is the model with lost sales. Consider the policy $\zeta(x) := 0$ that never orders. Then
        \begin{align*}
            v_\alpha^0(0) \leq \EE_0^\zeta \sum_{t=0}^{\infty} \alpha^t h(0) = (1-\alpha)^{-1} h(0) < +\infty,
        \end{align*}
        so $0 \in \dom v_\alpha^0$; hence, $[0, \tilde{x}] \in \dom v_\alpha^0$ by convexity.
    \end{enumcase}

    Thus, in all models we have shown that $\XX_\ell \subset \dom v_\alpha^0$.
    We next show that $\XX_r \subset \dom v_\alpha^0$. Again, there are two cases.
    \begin{enumcase}
        \item $\XX_r = [\tilde{x}, +\infty)$. This includes the models (U, BO). Let $x_r := \sup \dom v_\alpha^0$, and suppose $x_r < +\infty$. Then, starting at each $x \geq x_r$ under an optimal policy, the inventory level will never decrease past $x_r$; hence $\EE v_\alpha^0(x_r - D) = +\infty$ implies that $P(D=0) = 1$. Then $v_\alpha^{0}(x) \leq h(x) (1-\alpha)^{-1} < +\infty$, a contradiction. Thus $[\tilde{x}, +\infty) \subset \dom v_\alpha^0$.
        \item $\XX_r = [\tilde{x}, \bar x]$. This includes the models (BS, BOS). We estimate $v_\alpha^0(\bar x)$. Define the random variable $N = \inf\set{t \in \NN : x_t  \in \XX_{\ell}}$. Since $0 < \EE D < +\infty$, it follows from renewal theory (see, e.g.,~\cite[Proposition 3.2.2]{ross_1996_stochastic}) that $\EE N < +\infty$. Therefore, define the nonstationary policy $\phi$, which orders nothing until the time the inventory level first crosses into $\XX_\ell$, and then reverts to an optimal policy $\phi^*$. Then,
        \begin{equation*}
            v_\alpha^0(\bar x) \leq \EE_{\bar x}^\phi \sum_{t=0}^{N-1} \alpha^t h(x_{t+1}) + \alpha^{N} \EE_0^{\phi^*} v_0^{\alpha}(x_{N}) \leq \frac{1}{1-\alpha} h(\bar x) +  \alpha^{N} \EE_{\bar x}^{\phi^*} v_0^{\alpha}(x_{N}) < +\infty,
        \end{equation*}
        so $\bar x \in \dom v_\alpha^0$. Thus, by convexity we obtain $[\tilde{x}, \bar x] \subset \dom v_\alpha^0$.
    \end{enumcase}

    Thus, in all models we have shown that $\XX_r \subset \dom v_\alpha^0$. Therefore, $\dom v_\alpha^0 = \XX$, and continuity follows convexity, lower semicontinuity, and finiteness of $v_\alpha^0$.
  \end{proof}

\begin{proof}[Proof of Theorem~\ref{thm:continuity_of_infinite_horizon_values}]
  Let us show that $v_\alpha$ is finite for the models (U, BS). For the model (BS) the maximum storage size is $\bar x > 0$. The backorder rule is defined as $T(x) := 0 \lor x$ for the model with lost sales and $T(x) := x$ for the model with backorders. Hence, in both models the policy $\phi(x) := -x \lor 0$ is feasible. (In the model with lost sales, $\phi$ is simply the policy that never orders.) For any sequence of inventory levels $\set{x_t}_{t=0}^{\infty}$ under the policy $\phi$, we observe three inequalities: (i) $\phi(x_{t+1}) \leq D_{t+1}$, because the order size never exceeds the demand of the previous stage; (ii) $0 \leq x_t + \phi(x_t) \leq (0 \lor x_0)$ for each $t=0,1,\dots$, because $\phi$ never orders at positive inventory levels; and (iii) $h(T(x_t + \phi(x_t) - D_{t+1})) \leq h(x_t + \phi(x_t) - D_{t+1})$, since $h(T(x)) \leq h(x)$ for any inventory level $x \in \RR$. These imply that the estimate
  \begin{align*}
    v_\alpha^\phi(x)
    &\leq \EE_x^\phi \sum_{t=0}^{\infty} \alpha^t \left[K + \bar c D_{t+1} + h(T(x_t + \phi(x_t) - D_{t+1})) \right] \\%
    &\leq \EE_x^\phi \sum_{t=0}^{\infty} \left[K + \bar c D_{t+1} + h(-D_{t+1}) + h(x-D_{t+1}) \right] \\
    &\leq \frac{K+ \bar c\EE D + \EE [h(-D) + h(x-D)]}{1-\alpha},
  \end{align*}
  holds, for both forms of $T$, for each $x \in \XX$. Since the final term is finite, it follows that $v_\alpha(x) \leq v_\alpha^\phi(x) < +\infty$ for each $x \in \XX$ in the models (U, BS).

  The function $v_\alpha$ is the pointwise supremum of the finite horizon functions $v_{N, \alpha}$. From Theorem~\ref{thm:continuity_of_finite_horizon_values}, the continuity of each $v_{N, \alpha}$ implies that $v_\alpha$ is lower semicontinuous. To show that $v_\alpha$ is continuous, it therefore suffices to verify that $v_\alpha$ is upper semicontinuous. Let $v_\alpha^0$ denote the corresponding value function for the problem without setup costs as defined in Lemma~\ref{lem:no_setup_costs}. Then inequality~\eqref{eq:no_setup_inequality_infinite_horizon} implies that $v_\alpha$ is finite if and only if $v_\alpha^0$ is finite, which holds if and only if there is a point $x_0 \in \XX$ such that $v_\alpha^0(x_0) < +\infty$. Thus $v_\alpha$ is either finite on $\XX$, or $v_\alpha(x) \equiv +\infty$ on $\XX$.

  Define $G_\alpha^0(x) := \bar c x + \EE h(T(x-D)) + \alpha \EE v_\alpha^0(T(x-D))$, and we observe that $G_\alpha^0$ is convex and lower semicontinuous. The inequality~\eqref{eq:no_setup_inequality_infinite_horizon} implies that $G_\alpha^0(x) \leq G_\alpha(x) \leq G_\alpha^0(x) + \alpha(1-\alpha)^{-1} K$.  Therefore, $G_\alpha$ is finite if and only if $G_\alpha^0$ is finite. Suppose $v_\alpha^0$ is finite. We show that this implies that $G_\alpha^0$ is finite. Because $G_\alpha$ is finite for all $x \in \XX$ in the models (U, BS), there is no loss in generality to assume that $\bar a < +\infty$. Suppose $\tilde{x}\in \XX$ such that $G_\alpha(\tilde{x}) = + \infty$. Since $v_\alpha(\tilde{x}) < +\infty$, there is an $a > 0$ such that $G_\alpha(\tilde{x}+ a) < +\infty$. Since $\dom G_\alpha^0$ is a convex set, this implies that $\tilde{x} \leq \inf \dom G_\alpha^0$. On the other hand, let $x < \tilde{x} - \bar a$. Since $v_\alpha^0(x) < +\infty$, it follows that there is an action $a \in [0, \bar a]$ such that $G_\alpha(x + a) < +\infty$ and $x + a < \tilde{x}$. But this violates the convexity of $\dom G_\alpha^0$. Hence, $G_\alpha^0$ is finite everywhere. Conversely, if $G_\alpha^0 \equiv +\infty$, then evidently $v_\alpha^0 \equiv +\infty$. This implies, in view of the definition of $G_\alpha^0$ that $v_\alpha^0$ is finite if and only if the function $x \mapsto \EE v_\alpha^0(T(x-D))$ is finite. Similarly, the function $v_\alpha$ is finite if and only if the function $x \mapsto \EE v_\alpha(T(x-D))$ is finite. If $v_\alpha$ is infinite, then it is continuous. For the remaining proof, we shall assume instead that $v_\alpha$ is finite.

  For $n=0, 1, \dots$, we define the functions $f_n : \XX \to \RR$ by the  equations
  \begin{align*}
    f_0(x) &= v_\alpha^0(x) + (1-\alpha)^{-1}K , &&\qquad x \in \XX, \\
    f_{n+1}(x) &= \min_{a \in A(x)} \set{c(x,a) + \alpha \EE f_n(T(x+a-D))}, &&\qquad x \in \XX,
  \end{align*}
  and, in view of the convergence of value iterations, we observe that $v_\alpha(x) = \inf_{n =0,1,\dots} f_n(x)$ for each $x \in \XX$. In fact, $f_0$ satisfies the property that $f_{n} \leq f_0$ for each $n = 0, 1, 2, \dots$. Indeed, if $f_n \leq f_0$ for some $n$, then
  \begin{align*}
    f_{n+1}(x) &= \min_{a \in A(x)} \set{c(x,a) + \alpha \EE f_n(T(x+a-D))} \\
    &\leq \min_{a \in A(x)} \set{c(x,a) + \alpha \EE f_0(T(x+a-D))} \\
    &\leq \min_{a \in A(x)} \set{K+c^0(x,a) + \alpha \EE f_0(T(x+a-D))} \\
    &= f_0(x),
  \end{align*}
  where $c^0(x,a)$ defined in~\eqref{eq:no_setup_cost}. For convenience, denote $F_n(x) := \EE f_n(T(x-D))$. If each $f_n$ is continuous, then $v_\alpha$ is the pointwise infimum of continuous functions; hence, upper semicontinuous. As such, we shall show that $f_n$ is continuous. Since $v_\alpha^0$ and is continuous, it follows that $f_0$ is continuous. Furthermore, since $v_\alpha^0$ is convex and finite, the function $F_0(x)$ is convex and finite; hence continuous.

  Suppose $f_n$ is continuous. We show that $F_n$ is continuous. Since $f_n \leq f_0 < +\infty$ for each $n$, it follows that $F_n \leq F_0 < +\infty$ as well. Let $x_k \to x_0$ be a convergent sequence in $\XX$. Since $F_0$ is finite and continuous, it follows that $F_0(x_k) \to F_0(x)$, and the dominated convergence theorem implies that $F_n(x_k) \to F_n(x)$. Thus $F_n$ is continuous. Denote $c_n(x,a) := c(x, a) + \alpha F_n(x+a)$. Since $F_n$ is continuous, the function $c_n$ is lower semicontinuous. Further, since $v_\alpha \leq f_n$ for each $n = 0, 1, \dots$, and since~\eqref{eq:g_right_limit} implies $v_\alpha$ is inf-compact, it follows that $f_n$ is inf-compact; hence, $c_n$ is $\KK$-inf-compact.

  Let $x_0 \in \XX$ and $a_0 \in A(x_0)$ be arbitrary. Consider the policy $\phi(x) := a_0 \lor ((x-x_0) \land 0)$, which has the properties:
  \begin{inparaenum}[(i)]
    \item $\phi(x) \in A(x)$ for each $x \in \XX$,
    \item $\phi$ is continuous as a function $\XX \to \AA$,
    \item $\phi(x) \in A(x)$ for each $x \in \XX$, and
    \item $\phi(x_0) = a_0$.
  \end{inparaenum}
  We claim that $c_n(x_k, \phi(x_k)) \to c_n(x_0, a_0)$. Indeed, if $a_0 = 0$, then $\phi(x_k) = 0$ for each $k=1,2,\dots$, which implies that
  \begin{equation*}
    \lim_{k \to \infty} c_n(x_k, \phi(x_k)) = \lim_{k \to \infty} \EE h(T(x_k - D)) + \alpha F_n(x_k) = \EE h(T(x_0 - D)) + \alpha F_n(x_0) = c_n(x_0, a_0).
  \end{equation*}
  Alternatively, if $a_0 > 0$, then the sequence $\phi(x_k)$ is eventually positive, which implies that
  \begin{align*}
    \lim_{k \to \infty} c_n(x_k, \phi(x_k)) &= \lim_{k \to \infty} K + \bar c \phi(x_k) + \EE h(T(x_k + \phi(x_k) - D)) + \alpha F_n(x_k + \phi(x_k)) \\
    &= K + \bar c \phi(x_k) + \EE h(T(x_0 + \phi(x_k) - D)) + \alpha F_n(x_0 + \phi(x_k)) = c_n(x_0, a_0).
  \end{align*}
  In either case, we find that $c_n(x_k, \phi(x_k)) \to c_n(x_0, a_0)$. According to~\cite[Theorem 12]{feinberg_2021_continuity}, the function $c_n$ is feasible path transfer upper semicontinuous, and according to~\citet[Theorem 7]{feinberg_2021_continuity}, the function $f_{n+1}$ is continuous. Therefore, $v_\alpha$ is upper semicontinuous; hence continuous as desired.

  We now show that $G_\alpha$ is continuous. For $x_0 \in \XX$, we define the function $g_\alpha(x) = v_\alpha(x \lor (x_0 + 1)) + \bar c (x \lor (x_0 + 1))$, which is continuous and bounded. Therefore, $x \mapsto \EE g_\alpha(T(x-D))$ is continuous. Furthermore, $G_\alpha(x) = (1-\alpha) \bar c x + \EE h(T(x-D)) + \alpha \EE g_\alpha(T(x-D)) + \alpha \bar c \EE [D]$ is a sum of continuous functions for each $x  \in (-\infty, x_0+1]$, so $G_{\alpha}$ is continuous on $(-\infty, x_0+1]$, and in particular $G_\alpha$ is continuous at $x_0$. Since $x_0 \in \XX$ was arbitrary, $G_\alpha$ is continuous on $\XX$. Finally, we observe that $\EE v_\alpha(T(x-D)) = \alpha^{-1}(G_\alpha(x) - \bar c x - \EE h(T(x-D)))$ is the sum of continuous functions, which implies that $x \mapsto \EE v_{\alpha}(T(x-D))$ is continuous.
\end{proof}
}

The following example demonstrates that it is possible for the inventory models (BO, BOS) to have $v_\alpha(x) = +\infty$ for each $x \in \XX$.

\begin{example}
	Consider the models (BO, BOS) with $D$ is almost surely constant $d>\bar{a}$, and with the holding/backorders cost function $h(x) = \alpha^{-x^2}$. For any feasible policy $\phi$,
	\begin{equation}
		v^\phi_\alpha(x)
		= \sum_{t=0}^{\infty} \alpha^t c(x_t, a_t)
		\geq \sum_{t=0}^{\infty} \alpha^t h(x_{t+1})
		= \sum_{t=0}^{\infty} \alpha^{t-x_{t+1}^2}.
		\label{eq:deterministic_divergent_series}
	\end{equation}
	For $t=0,1,\dots$, it follows that
	\begin{equation}
		x-dt \leq x_t - d \leq x_{t+1} \leq x_t - (d-\bar{a}) \leq x - (d-\bar{a})t,
	\end{equation}
	so in the limit
	\begin{equation}
		\lim_{t \to \infty} t - x_{t+1}^2 \leq \lim_{t \to \infty} t -[x - (d-\bar{a})t]^2 = -\infty,
	\end{equation}
	\revision{which implies that $\alpha^{t-x_{t+1}^2} \to +\infty$. This means that $v_\alpha^\phi(x) = +\infty$ for each $x \in \XX$ and for each policy $\phi$.} Hence, $v_\alpha(x) = +\infty$ for each $x \in \XX$.
  \label{ex:infinite_values}
\end{example}
\revision{
  \section{Proof of Theorem~\ref{thm:bounded_storage_s_S_policies}}
    \label{sec:bounded_storage}
}
In this section we establish the correctness of Theorem~\ref{thm:bounded_storage_s_S_policies}. We consider the models (U) with unbounded order sizes and storage capacity and (BS) with unbounded order sizes and bounded storage capacity.  Following the presentation of~\citet[Chapter 4]{rockafellar_2009_variational}, denote the collection of all infinite subsets of $\NN$ by $\mathcal{N}^{\#}_{\infty} = \set{N \subset \NN : \text{$N$ is infinite}}$, and let $X$ be a subset of a Euclidean space. For a sequence of sets $C_n \subset X$, the outer limit of $\set{C_n}_{n \in \NN}$ is defined as
\begin{align}
  \limsup_{n \to \infty} C_n &= \set{x \in X : \exists N \in \mathcal{N}_{\infty}^{\#}, \exists x_n \in C_n (n \in N), x_n \to x},
  \label{eq:pk_limsup}
\end{align}
which is composed of all limit points of sequences $\set{x_n}_{n \in \NN}$ with $x_n \in C_n$ for each $n \in \NN$. We observe that the outer limit is always a closed subset of $X$.

The following proposition is a classic fact about $K$-convex functions that connects them to $(s,S)$ policies; see, e.g.,~\citet[Lemma 4.2.1(d)]{bertseks_2005_dynamic} or~\citet[Lemma 9.3.2(d)]{simchi_levi_2014_logic}

\begin{proposition}[{\citet[Lemma 4.2.1(d)]{bertseks_2005_dynamic}}]
  Let $f : \RR \to \RR$ be a continuous, inf-compact $K$-convex function. Let $S$ and $s$ be defined as in equations~\eqref{eq:big_S} and~\eqref{eq:little_s}, respectively.
  Then
  \begin{enumerate}
    \item $f(S) + K < f(x)$ for all $x < s$, \label{item:order_to_S_i}
    \item $f(x)$ is decreasing on $(-\infty, s]$, so $f(s) < f(x)$ for all $x < s$, \label{item:order_to_S_ii}
    \item $f(x) \leq f(z) + K$ for all $s \leq x \leq z$. \label{item:order_nothing}
  \end{enumerate}
  \label{prop:s_S_policies}
\end{proposition}
Proposition~\ref{prop:s_S_policies} is the mechanism by which $(s,S)$ polices can be shown to be optimal for $K$-convex functions $f$. In particular, Statements~\ref{item:order_to_S_i} and~\ref{item:order_to_S_ii} imply that for each $x < s_f$, it is preferable to order $S_f-x$ than to order nothing, and Statement~\ref{item:order_nothing} implies that for each $x \geq s_f$, it is always optimal not to order.

A an important property of $K$-convex functions is that the function
\begin{equation}
  g(x) := \min_{a \geq 0} \set{K \one{a > 0} + f(x + a)}
  \label{eq:minimum_k_convex_unbounded_orders_unbounded_storage}
\end{equation}
is $K$-convex, if $f$ is $K$-convex; see, e.g.,~\citet[Proposition 9.3.3]{simchi_levi_2014_logic}. The immediate application of this fact is to the unbounded orders, unbounded storage inventory model (U), but~\eqref{eq:minimum_k_convex_unbounded_orders_unbounded_storage} also holds when the minimization is taken over all $a \in [0, 0 \lor (\bar x - x)]$, which is shown in the following lemma.

\begin{lemma}
  Let $f : \RR \to \RR$ be a $K$-convex function. Then the function
  \begin{equation}
    g(x) := \min_{0 \leq a \leq \bar x - x} \set*{K \one{a > 0} + f(x+a)}
    \label{eq:minimum_k_convex_unbounded_orders_bounded_storage}
  \end{equation}
  is also $K$-convex.
  \label{lem:minimum_k_convex}
\end{lemma}

\begin{remark}
  In the formulation of Proposition 9.3.3 in~\citet{simchi_levi_2014_logic}, the coefficient in front of the indicator $\one{a > 0}$ is an arbitrary $Q>0$ instead of $K.$ The similar generalization  holds for \eqref{eq:minimum_k_convex_unbounded_orders_bounded_storage}, but it is not used in the current paper.
\end{remark}

\begin{proof}[Proof of Lemma~\ref{lem:minimum_k_convex}]
  The proof of Proposition 9.3.3 in~\citet{simchi_levi_2014_logic} actually suffices for the modified equation~\eqref{eq:minimum_k_convex_unbounded_orders_bounded_storage}, but we reproduce it in order to deal with the added constraint explicitly. The property that the proof in~\citet[Proposition 9.3.3]{simchi_levi_2014_logic} requires throughout the argument is that, if $x_0 \leq x_1$ and $a \geq 0$, then $a + (x_1 - x_0) \geq 0$. The related property holds for~\eqref{eq:minimum_k_convex_unbounded_orders_bounded_storage}, since if $0 \leq a \leq \bar x - x_1$, then $0 \leq a + (x_1 - x_0) \leq \bar x - x_0$. This leaves the basic argument unchanged. Let $E = \set{x : g(x) = f(x)}$ and $O = \set{x : g(x) < f(x)}$. Let $x_0 \leq x_1$, let $\theta \in [0,1]$, and denote $x_\theta := (1-\theta) x_0 + \theta x_1$. We consider four cases.
  \begin{enumcase}
    \item $x_0, x_1 \in E$. Then since $a=0$ is feasible for all $x$,
    \begin{equation*}
      g(x_\theta)
      \leq f(x_\theta)
      \leq (1-\theta) f(x_0) + \theta f(x_1) + \theta K
      = (1-\theta) g(x_0) + \theta g(x_1) + \theta K.
    \end{equation*}

    \item $x_0, x_1 \in O$. Fix $0 \leq a_0 \leq \bar x - x_0$ and $0 \leq a_1 \bar x - x_1$ such that $g(x_0) = K + f(x_0+a_0)$ and $g(x_1) = K + f(x_1+a_1)$. Now, the inequality $x_0 + a_0 \leq x_1 + a_1$ holds. Indeed, if $x_0 + a_0 < x_1$, this is immediate. Otherwise, if $x_0 + a_0 \geq x_1$, it follows that $a_0 + (x_0 - x_1) \leq \bar x - x_1$. Therefore, $a_1 \geq x_0 + a_0 - x_1$, and the inequality holds. Furthermore, $x_\theta \leq x_\theta + a_\theta$. Thus
    \begin{align*}
      g(x_\theta)
      &\leq K + f(x_\theta + a_\theta)\\
      &\leq (1-\theta) (K + f(x_0+a_0)) + \theta (K + f(x_1+a_1)) + \theta K\\
      &= (1-\theta) g(x_0) + \theta g(x_1) + \theta K.
    \end{align*}

    \item $x_0 \in E$, $x_1 \in O$. Fix $0 \leq a_1 \leq \bar x - x_1$ such that $g(x_1) = K + f(x_1+a_1)$, and fix $\nu$ such that $x_\theta = (1-\nu) x_0 + \nu (x_1 + a_1)$. We observe that $\nu \leq \theta$, and
    \begin{align*}
      g(x_\theta)
      &\leq f(x_\theta) \\
      &\leq (1-\nu) f(x_0) + \nu f(x_1+a_1) + \nu K \\
      &= (1-\theta) g(x_0) + \theta g(x_1) + (\theta - \nu) (f(x_0) - f(x_1+a_1) - K) \\
      &\leq (1-\theta) g(x_0) + \theta g(x_1) \\
      &\leq (1-\theta) g(x_0) + \theta g(x_1) + \theta K,
    \end{align*}
    where the third inequality follows from the following observation: since $f(x_0) = g(x_0) \leq K + f(x_0 + a)$ for all $0 \leq a \leq \bar x - x_0$, and since $0 \leq a_1 + (x_1-x_0) \leq \bar x - x_0$, then $f(x_0) \leq K + f(x_1 + a_1)$.

    \item $x_0 \in O$, $x_1 \in E$. Fix $0 \leq a_0 \leq \bar x - x_0$ such that $g(x_0) = K + f(x_0+a_0)$. We observe that $0 \leq x_1 - x_0 \leq \bar x - x_0$. Therefore, $g(x_0) \leq K + f(x_1)$, which implies $f(x_0 + a_0) \leq f(x_1)$. Now, if $x_\theta \leq x_0 + a_0$, then $0 \leq a_0 + x_0 - x_{\theta} \leq \bar x - x_{\theta}$, which implies
    \begin{align*}
      g(x_\theta)
      &\leq K + f(x_0+a_0) \\
      &= (1-\theta) (K + f(x_0+a_0)) + \theta f(x_1) + \theta(K + f(x_0) - f(x_1)) \\
      &\leq (1-\theta) g(x_0) + \theta g(x_1) + \theta K,
    \end{align*}
    where the final inequality follows from $f(x_0 + a_0) \leq f(x_1)$.	On the other hand, if $x_\theta \geq y_0$, then fix $\nu \in [0,1]$ such that $x_\theta = (1-\nu) y_0 + \nu x_1$, and we observe that $\nu \leq \theta$. Then it follows that
    \begin{align*}
      g(x_\theta)
      &\leq f(x_\theta) \\
      &\leq (1-\nu) f(y_0) + \nu f(x_1) + \nu K \\
      &= (1-\theta) g(x_0) + \theta g(x_1) + \nu K + (\theta - \nu) (f(y_0) - f(x_1)) - (1-\theta) K \\
      &\leq (1-\theta) g(x_0) + \theta g(x_1) + \theta K,
    \end{align*}
    where the final inequality follows again from $f(x_0 + a_0) \leq f(x_1)$.
  \end{enumcase}

  The above four cases exhaust all combinations of $x_0$ and $x_1$, so $g$ is indeed $K$-convex.
\end{proof}

The $K$-convexity of $g$ defined in~\eqref{eq:minimum_k_convex_unbounded_orders_unbounded_storage} has direct applications to the structure of optimal policies for the unbounded orders, unbounded storage inventory model (U) with backorders. Lemma~\ref{lem:minimum_k_convex} now extends these applications to the unbounded orders, bounded storage model (BS) with backorders. Throughout the remainder of this section,the constants $\alpha^\star$ and $N_\alpha$ defined in~\eqref{eq:alpha_star} and~\eqref{eq:N_alpha} will be used.

\begin{lemma}
  For the models (U, BS) with backorders, the functions $G_{t,\alpha},$ $t=0,1,\ldots,$ and $G_{\alpha}$, $\alpha\in [0,1),$ are $K$-convex,  and
  \begin{enumerate}
    \item if $\alpha \leq \alpha^\star$, then $G_{t,\alpha}$ and $G_{\alpha}$ are convex and nondecreasing functions for $t=0,1,\ldots;$
    \item if $\alpha > \alpha^\star$, then $G_{t,\alpha}$ is convex and nondecreasing for $t=0, 1, \dots, N_\alpha$, and $G_{t,\alpha}$  and $G_{\alpha}$ inf-compact for $t=N_\alpha+1, N_\alpha+2, \dots$
    \label{item:eventually_k_convex_ii}
  \end{enumerate}
  \label{lem:eventually_k_convex}
\end{lemma}

\begin{proof}
  Let $\alpha \in [0,1)$ be arbitrary.
  Denote by $v_{t,\alpha}^{\mathrm{U}}, v_{\alpha}^{\mathrm{U}}, G_{t,\alpha}^{\mathrm{U}}, G_{\alpha}^{\mathrm{U}}$ the functions defined in~\eqref{eq:finite_horizon_bellman}-\eqref{eq:infinite_horizon_g} for the model (U), and similarly $v_{t,\alpha}^{\mathrm{BS}}, v_{\alpha}^{\mathrm{BS}}, G_{t,\alpha}^{\mathrm{BS}}, G_{\alpha}^{\mathrm{BS}}$ for the model (BS). Then the inequalities
  \begin{equation}
    v_{t,\alpha}^{\mathrm{BS}}  \geq v_{t,\alpha}^{\mathrm{U}}, \quad
    v_{\alpha}^{\mathrm{BS}}    \geq v_{\alpha}^{\mathrm{U}}, \quad
    G_{t, \alpha}^{\mathrm{BS}} \geq G_{t,\alpha}^{\mathrm{U}}, \quad
    G_{\alpha}^{\mathrm{BS}}    \geq G_{\alpha}^{\mathrm{U}}
  \end{equation}
  hold for each $t = 0, 1, \dots$, \revision{since the feasible orders for the model (BS) are always feasible for the model (U). If $G_{t,\alpha}^{\mathrm{U}}$ is convex and nondecreasing for some $\alpha$ and $t$, then according to~\citet[Lemma 4.6]{feinberg_2017_structure}, then the optimality equation~\eqref{eq:finite_horizon_bellman} is achieved with $a=0$, and so $G_{t,\alpha}^{\mathrm{BS}} = G_{t,\alpha}^{\mathrm{U}}$ and $G_{\alpha}^{\mathrm{BS}} = G_{\alpha}^{\mathrm{U}}$, and hence each of these functions is convex (and $K$-convex).  From~\eqref{eq:alpha_star}, it follows that $\alpha^\star < 1$. We consider two cases.
  \begin{enumcase}
    \item $\alpha^\star \geq 0$. Then according to~\citet[Theorem 4.2(ii)]{feinberg_2017_structure}, the function $G_{N_\alpha, \alpha}^{\mathrm{U}}$ is convex (hence, $K$-convex) and nondecreasing for $t=0, 1, \dots, N_\alpha-1$, and $G_{t,\alpha}^{\mathrm{U}}$ and $G_{\alpha}^{\mathrm{U}}$ are convex (hence, $K$-convex), and nondecreasing for $t= N_\alpha, N_\alpha+1, \dots$.
    \item $\alpha^\star < 0$. Then according to~\citet[Lemma 6.11]{feinberg_2018_convergence}, the function $G_{0, \alpha}^{\mathrm{U}}$ is $K$-convex and inf-compact for each $\alpha \geq 0$.
  \end{enumcase}

  In either case, if $\alpha \leq \alpha^\star$, then $G_{t,\alpha}^{\mathrm{U}}$ is convex and nondecreasing for each $t = 0, 1, \dots$, which implies that the optimality equation~\eqref{eq:finite_horizon_g} is satisfied by the order $a=0$ for each $x \in \XX$. Since $a=0$ is feasible in the model (BS), the equality $G_{t,\alpha}^{\mathrm{U}} = G_{t,\alpha}^{\mathrm{BS}}$ holds, and $G_{t,\alpha}^{\mathrm{BS}}$ is convex and nondecreasing for each $t = 0, 1, \dots$. The similar analysis holds for $G_{\alpha}^{\mathrm{U}}$ and $G_{\alpha}^{\mathrm{BS}}$. On the other hand, if $\alpha^\star < \alpha$, then according to~\citet[proof of Theorem 4.2]{feinberg_2017_structure} for each $t = N_{\alpha}+1,N_{\alpha}+2, \dots$ the functions $G_{t,\alpha}^{\mathrm{U}}$ are $K$-convex such that $\lim_{\abs{x} \to +\infty} G_{t,\alpha}(x) = +\infty$. Since $G_{t,\alpha}^{\mathrm{BS}} \geq G_{t,\alpha}^{\mathrm{U}}$, it follows that $G_{t,\alpha}^{\mathrm{BS}}$ is $K$-convex and $\lim_{x \to -\infty} G_{t,\alpha}^{\mathrm{BS}}(x) = +\infty$. Therefore, $G_{\alpha}^{\mathrm{BS}}$ is also $K$-convex.}
\end{proof}

\revision{%
\begin{lemma}
  Suppose $0 \leq \alpha^\star < \alpha$. Then the equality~\eqref{eq:N_alpha_formula}
  holds, and $1 \leq N_\alpha < +\infty$.
  \label{lem:N_alpha_equation}
\end{lemma}
\begin{proof}
  Let $k_h := -\lim_{x \to -\infty} h(x) / x$, so that $k_h = \bar c (1-\alpha^\star)$. Then the statement $\alpha^\star \geq 0$ is equivalent to the fact that $k_h \leq \bar c$. There is no loss o generality to assume that $h(0) = 0$. Then it follows from convexity that $h(x) \leq k_h \abs{x}$ for each $x < 0$. From the definition of $f_{t,\alpha}$ in~\eqref{eq:f_t_alpha}, there exists $t \in \NN$ such that $\lim_{x \to -\infty} f_{t,\alpha}(x) = +\infty$, if and only if $\lim_{x \to -\infty} f_{t,\alpha}(x)/x < 0$.

  The inequality
  \begin{align*}
      \lim_{x \to -\infty} \frac{f_{t,\alpha}(x)}{x} = \bar c + \lim_{x \to -\infty}\sum_{i=0}^{t} \alpha^i \frac{\EE h(x-\mathbf{S}_{i+1})}{x}
      \geq \bar c +  \lim_{x \to -\infty} \sum_{i=0}^{t} \alpha^i \frac{h(x-\EE D)}{x}
      = \bar c -k_h \sum_{i=0}^{t} \alpha^i
  \end{align*}
  follows from Jensen's inequality and from $h(x-\mathbf{S}_{i+1}) \geq h(x-D)$ for all $x < 0$ and $i=0,1,\dots$. Conversely, the inequality
  \begin{align*}
   \lim_{x \to -\infty} \frac{f_{t,\alpha}(x)}{x} \leq \bar c + \lim_{x \to -\infty} \sum_{i=0}^{t} \alpha^i \frac{k_h(x-\EE\mathbf{S}_{i+1})}{x}
   = \bar c -k_h \sum_{i=0}^{t} \alpha^i
  \end{align*}
  follows from the fact that $h(x) \leq -k_h x$ for all $x < 0$.  Therefore, we obtain the equality
  \begin{equation*}
      \lim_{x \to -\infty} \frac{f_{t,\alpha}(x)}{x} = \bar c -k_h \sum_{i=0}^{t} \alpha^i.
  \end{equation*}
  Furthermore, by taking the limit $t \to \infty$ we obtain the inequality
  \begin{equation*}
      \bar c - \frac{k_h}{1-\alpha} = \bar c\left(1 - \frac{1-\alpha^\star}{1-\alpha}\right) < 0,
  \end{equation*}
  which follows from $\alpha^\star < \alpha$.
  As such, there exists $t \in \NN$ for which the finite inequality holds as well. Since $N_\alpha \leq t$ if and only if (after some algebra) $\frac{\alpha^\star}{1-\alpha^\star} < \sum_{i=1}^{t} \alpha^i$, we thus find that $N_\alpha$ is the minimum of all such $t$.
\end{proof}
}

\begin{proof}[Proof of Theorem~\ref{thm:bounded_storage_s_S_policies}]
  \revision{%
    We first consider $0 \leq \alpha \leq \alpha^\star$. According to Lemma~\ref{lem:eventually_k_convex}, the functions $G_{t,\alpha}$ and $G_\alpha$ are convex and nondecreasing. Therefore,
    \begin{align*}
      v_{t, \alpha}(x) &= \min_{a \in A(x)} \set{K \one{a > 0} + G_{t, \alpha}(x+a)} -\bar c x = G_{t,\alpha}(x)  -\bar c x,  \\
      v_{\alpha}(x) &=  \min_{a \in A(x)} \set{K \one{a > 0} + G_{\alpha}(x+a)} -\bar c x = G_{\alpha}(x) -\bar c x,
    \end{align*}
    which means that for each $x \in \XX$, the optimal action is $a=0$. Hence the no-order policy is $N$-stage optimal and discount-optimal.

    We now consider $\alpha > \alpha^\star$. Lemma~\ref{lem:N_alpha_equation} establishes Statement~\ref{item:N_alpha_equation}.
  }
  For the model (U), the theorem follows from~\citet[Theorem 4.3]{feinberg_2017_structure}, so we establish the result for the model (BS).                                                                                                                                                                                                                     t                                           t
  To prove Statement~\ref{item:discounted_finite_s_S}, we apply Lemma~\ref{lem:eventually_k_convex} to conclude that $G_{N, \alpha}$ is $K$-convex and inf-compact for each $\alpha \in [\alpha^\star, 1)$ and $N = N_\alpha+1, N_\alpha+2, \dots$. Then by Proposition~\ref{prop:s_S_policies}, for each $t = 0, 1, \dots, N-N_\alpha-1$, we set $s_{t,\alpha} := s_{G_{N-t-1,\alpha}}$ and $S_{t,\alpha} := S_{G_{N-t-1, \alpha}}$, since each $G_{N_\alpha+1, \alpha}, \dots, G_{N-1, \alpha}$ is inf-compact and $K$-convex. For $t=N-N_\alpha-1, \dots, N-1$, Lemma~\ref{lem:eventually_k_convex} implies that the function $G_{t,\alpha}$ is convex and nondecreasing, so the policy which never orders is optimal.

  For $\alpha^\star$ defined in~\eqref{eq:alpha_star}, Lemma~\ref{lem:eventually_k_convex} implies that $G_\alpha$ is $K$-convex and inf-compact. Therefore, by Proposition~\ref{prop:s_S_policies}, we set $s_\alpha := s_{G_\alpha}$ and $S_\alpha := S_{G_\alpha}$, and the stationary policy $(s_\alpha, S_\alpha)$ is optimal. Furthermore, since each $G_{t,\alpha}$ is continuous with $G_{t,\alpha} \uparrow G_{\alpha}$, according to~\citet[Proposition 7.4(d), Theorem 7.33]{rockafellar_2009_variational}, the inclusion
  \begin{equation}
    \set{(s_{N,\alpha}, S_{N,\alpha})}_{N=N_\alpha+1}^{\infty} \subset \limsup_{N \to \infty} \set{x \in \XX : G_{N,\alpha}(x) \leq K + \inf G_{N,\alpha}} \subset \set{x \in \XX : G_\alpha(x) \leq K + \inf G_{\alpha}}
    \label{eq:bounded_s_S_policies}
  \end{equation}
  holds. The set on the right-hand side of~\eqref{eq:bounded_s_S_policies} is compact, so the sequence on the left-hand side of~\eqref{eq:bounded_s_S_policies} is contained in a compact set; hence, every limit point of $\set{(s_{N,\alpha}, S_{N,\alpha})}_{t=N_\alpha+1}^{\infty}$ is discount optimal, and Statement~\ref{item:discounted_infinite_s_S} is proved.
\end{proof}

\bibliography{ArxivFK_InventoryControl_20220725_2352}
\end{document}